\documentclass{article}%
\usepackage{amsmath}
\usepackage{amsfonts}
\usepackage{amssymb}
\usepackage{graphicx}%
\begin{document}

\title{Perturbed Operators on Banach Spaces}
\author{Jos\'{e} Mar\'{\i}a Soriano Arbizu$^{\dag}$ and Manuel Ordo\~{n}ez Cabrera
\and Departamento de An\'{a}lisis Matem\'{a}tico, Facultad de Matem\'{a}ticas,
\and Universidad de de Sevilla, c/ Tarfia s/n, Sevilla \ 41012, Spain}
\date{}
\maketitle

\begin{abstract}
Let $X$ be a Banach space over $\mathbb{K}=\mathbb{R}$ or $\mathbb{C}$, and
let $f:=F+C$ be a weakly coercive operator from $X$ onto $X$, where $F$ is a
$C^{1}-$proper operator, and $C$ \ a $C^{1}-$compact operator. Sufficient
conditions are provided to assert that the perturbed operator $f$ is a
$C^{1}-$diffeomorphism. \ Three corollaries \ are given. The first one, when
$F$ is a linear homeomorphism. The second one, when $F$ is a $k-$contractive
perturbation of the identity. The third one, when $X$ is a Hilbert space \ and
$F$ a particular linear operator. The proof of our results is based on
properties of Fredholm operators, as well as on local and global inverse
mapping theorems, and the Banach fixed point theorem. \ As an \ application
two examples are given.

\end{abstract}

\footnotetext{{}The authors have been partially supported by the Plan Andaluz
de Investigaci\'{o}n de la Junta de Andaluc\'{\i}a FQM-127 Grant
P08-FQM-03543, and by MEC Grant MTM2015-65242-C2-1-P.
\par
\dag Corresponding author, E-mail: soriano@us.es}


\medskip

\begin{center}
\textbf{1. INTRODUCTION\ AND PRELIMINARIES}

\bigskip
\end{center}

The existence of $C^{r}-$diffeomorphisms, $0\leq r\leq\infty$, are used in
applications of functional analysis. For example, given an operator equation
$fu=v$ between Banach spaces $X,Y$, the knowledge that the operator
$f:X\rightarrow Y$ is a $C^{r}-$diffeomorphism naturally implies the existence
and uniqueness of solution $u$ in the initial space $X$ \ for any fixed $v$ of
the final space $Y$. In practice, operators can be perturbed\ by the addition
of other operators. The goal of our paper is to study the behavior of some
these perturbed operators. We show that when $Y=X$, and a $C^{1}-$proper
operator $F$ is perturbed \ by a $C^{1}-$compact operator $C$, the resultant
perturbed operator $f:=F+C$ is a global $C^{1}-$diffeomorphism, if some
regularity conditions are verified. Hence existence and uniqueness of solution
are verified \ for any fixed value $v$ of the final space. When the hypotheses
required to the operators are weakened but we suppose that a fixed $v$ is a
regular value for the perturbed operator, we only obtain that this particular
equation have at most finitely many solutions.

A first corollary of our main result shows that if $F$ is a linear
homeomorphism and $C$ a linear compact operator with $\left\Vert F\right\Vert
\neq\left\Vert C\right\Vert $ \ that $f$ is a global $C^{1}-$diffeomorphism. A
second corollary shows that a linear $k-$contractive \ plus a $C^{1}-$compact
perturbation of the identity is a global $C^{1}-$diffeomorphism, if some
regularity conditions are verified. A third corollary is given for $C^{1}%
-$compact perturbations of particular linear operators on Hilbert spaces.

The tools used are local and global inverse mapping theorems, and a corollary
of the Banach fixed point theorem [6], as well as properties of Fredholm
mappings [8].

We could have used degree theory [2] or continuation methods [1] trying to
prove the existence of finite number of solutions, but this does not guarantee
the desired uniqueness.

As an application of our results, two examples are furnished.

Other existence theorems have been given by the first author in finite
dimensional setting [4] and in infinite dimensional setting [5], and by both
authors in [3].\smallskip

\bigskip

\textbf{Definitions} [6-8].

The operator $T:M\subseteq X\rightarrow X$, where $\left(  X,d\right)  $ is a
metric space and $M$ is a closed nonempty subset of $X$, is called
$k-$\textit{contractive} iff the following holds: $d(Tx,Ty)\leq kd(x,y)$ for
all $x,y\in M$ and for a fixed $k$, $0\leq k<1$. We call to this constant $k$,
the contractive constant of the operator $T$.

In the following definitions, $X$ and $Y$ are Banach spaces over
$\mathbb{K}=\mathbb{R}$ or $\mathbb{C}$, if nothing else is specified.

Given two operators $f,g:X\rightarrow Y$, then $f$ is called a\textit{ compact
perturbation} of $g$ (a $k-$\textit{contractive perturbation}) provided that
there is a compact operator (a $k-$\textit{contractive operator}),
$h:X\rightarrow Y$ \ such that, $\ f=g+h$.

An operator $A:X\rightarrow Y$ is called \textit{weakly coercive} whenever
$\Vert A(u)\Vert$ $\rightarrow$ $\infty$ \ if and only if $\Vert
u\Vert\rightarrow\infty$.

The symbol $\dim$ means \textit{dimension}, $\operatorname*{codim}$ means
\textit{codimension},\textit{ }$\ker$ means \textit{kernel}, and \textit{
}$\operatorname{R}(L)$ stands for \textit{range} of the operator $L$.\ 

An operator $A:D(A)\subset X\rightarrow Y$ is said to be \textit{compac}t
whenever it is continuous and the image $A(B)$ is relatively compact (i.e. its
closure $\overline{A(B)}$ is compact in $Y$) for every bounded subset $B$
$\subset D(A)$. Obviously this property is equivalent to the following: If
$(u_{n})_{n\geq1}$ is a bounded sequence in $D(A)$, then there is a
subsequence $(u_{n_{k}})_{k\geq1}$ \ such that the sequence $(A(u_{n_{k}%
}))_{k\geq1}$ is convergent in $Y.$

The operator $A$ is said to be \textit{proper} whenever the pre-image
$A^{-1}(K)$ of every compact subset $K\subset Y$ is also a compact subset of
$D(A).$

The operator $A$ is said a \textit{submersion} at the point $u$ if it is a
$C^{1}-$mapping on a neighborhood of $u$, $A^{\prime}(u):X\rightarrow Y$ is
surjective, and $\ker(A%
\acute{}%
(u))$ \textit{\ }has topological complement in $X.$\ The point $u\in X$ is
called a \textit{regular point of }$A$ if and only if $A$ is a submersion at
$u$. The point $v\in Y$ is called a regular value if and only if $A^{-1}(v)$
is empty or consists solely of regular points.

That $L:X\rightarrow Y$ is a \textit{linear Fredholm} operator means that $L$
is linear and continuous and both numbers $\dim$($\ker(L)$) and
$\operatorname*{codim}(\operatorname*{R}(L))$ are finite, and therefore
$\ker(L)=:X_{1}$ is a Banach space and has topological complement $X_{2}$,
since $\dim(X_{1})$ is finite. The integer number $\operatorname{Ind}(L)$=
$\dim$($\ker(L$)) $-$ $\operatorname*{codim}(\operatorname*{R}(L))$ is called
the \textit{index} of $L$.

Let $A:D(A)\subseteq X\rightarrow$ $Y$. If $D(A)$ is open, then $A$ is said to
be a \textit{Fredholm} operator whenever both $A$ is a $C^{1}$-operator and
the \textit{Fr\'{e}chet derivative} of the operator $A$ at the point $x$,
$A^{\prime}(x):X\rightarrow Y$, is a Fredholm linear operator for all $x\in
D(A).$ If $\operatorname*{Ind}(A^{\prime}(x))$ is constant with respect \ to
\ $x\in D(A)$, then we call this number the index of $A$ and write it as
$\operatorname*{Ind}(A)$.

Two Banach spaces $X,Y$ are called \textit{isomorphic} if and only if there is
a linear homeomorphism (\textit{isomorphism}) $L:X\rightarrow Y$.

By $\mathcal{F}(X,Y)$, $\mathcal{L}(X,Y)$ and Isom$(X,Y)$ we denote,
respectively, the set of all linear Fredholm operators $A:X\rightarrow Y$, the
set of all linear continuous operators $L:X\rightarrow Y$, and the set of all
isomorphisms $L:X\rightarrow Y$.

Let $f:M\subseteq X\rightarrow Y$ with $M$ an arbitrary set. Then $f$ is
called a $C^{r}-$operator when it can be extended \ locally to a $C^{r}%
-$operator in the usual sense, wich means that, for each $x\in M$, there
exists an open neighborhood $U(x)$ such that $f$ can be extended to a $C^{r}%
-$operator on $U(x)$.

Let $M$ and $N$ be arbitrary sets in $X$ and $Y$ respectively. Let $r$ be
either a nonnegative integer or $\infty$. Then a mapping $\ f:M\rightarrow N$
is called a \ $C^{r}-$\textit{diffeomorphism} if and only if $f$ is bijective
and both $f$ and $\ f^{-1}$ are $C^{r}-$operators. A \textit{local} $C^{r}%
-$\textit{diffeomorphism }at $x_{0}$ is a $C^{r}-$\textit{diffeomorphism }from
some neighborhood $U(x_{0})$ in $X$ onto some neighborhood $U(f(x_{0}))$ in
$Y$.

For a map $r:U(0)\subset X\rightarrow Y$, we will write $r(x)=o(\left\Vert
x\right\Vert )$, $x\rightarrow0$, provided that $\frac{r(x)}{\left\Vert
x\right\Vert }\rightarrow0$ as $x\rightarrow0$.

By $\left(  \cdot\mid\cdot\right)  $ we denote the \textit{inner} product on a
Hilbert space X. By $\operatorname{span}(\Omega)$ we denote the Lebesgue
measure of the set $\Omega$. By $f_{y}(x,y)$ is denoted the partial derivative
of $f$ at $(x,y)$ with respect $y$.

\medskip

\textbf{Theorem 1}. (\textit{Global inverse mapping theorem of Banach and
Mazur}) [6]\textit{. Assume that }$0\leq r\leq\infty$. \textit{Let
\ }$f:X\rightarrow Y$ \textit{be a local} $C^{r}-$\textit{diffeomorphism},
\textit{at every point of }$X$. \textit{Then} $\ f$ \textit{is a} $C^{r}%
-$\textit{diffeomorphism if and only if} $f$ \textit{is proper}.\medskip

\textbf{Theorem 2}. (\textit{Local inverse mapping theorem}) [6]. \textit{Let}
$f:U(x_{0})\subseteq X\rightarrow Y$ \textit{be a} $C^{1}-$\textit{mapping.
Then} $f$ \textit{is a local }$C^{1}-$\textit{diffeomorphism at }$x_{0}%
$\textit{ if and only if} $f^{\prime}(x_{0}):X\rightarrow Y$ \textit{is
bijective}.\medskip

\textbf{Theorem 3}. (\textit{Banach's open mapping theorem}) [6]. \textit{If}
$A:X\rightarrow Y$ \textit{is linear, continuous,and surjective, then:} (i)
$A$ \textit{maps open sets into open sets}.(ii) \textit{If the inverse
operator }$A^{-1}:Y\rightarrow X$ \textit{exists, then it is continuous}%
.\medskip

\textbf{Theorem 4}.\textbf{ }[6].\textit{ Let }$g:D(g)\subset X\rightarrow
Y$\textit{ be a compact \ operator. Let }$a\in D(g).$\textit{ If the
derivative }$g^{\prime}(a)$\textit{ exists, then }$g^{\prime}(a)$\textit{ is
also a compact operator.\medskip}

\textbf{Theorem 5}. [8]. \textit{Let} $S$ $\in\mathcal{F}(X,Y)$. \textit{The
perturbed operator} $S+C$ \textit{\ satisfies} $S+C\in\mathcal{F}(X,Y)$
\textit{and }$\operatorname*{Ind}(S+C)=\operatorname*{Ind}(S)$ \textit{if}
$C\in\mathcal{L}(X,Y)$ \textit{and} $C$ \textit{is a compact \ operator}%
.\medskip

\textbf{Theorem 6}. [8]. \textit{Let} $S$ $\in\mathcal{F}(X,Y)$. \textit{Then
there exists a number} $\xi>0$ \textit{such that} $T\in\mathcal{F}(X,Y)$
\textit{and} $\operatorname*{Ind}(T)=\operatorname*{Ind}(S)$ \textit{for all
operators} $T\in\mathcal{L}(X,Y)$ \textit{with}$\left\Vert T-S\right\Vert
<\xi$.

\textbf{Theorem 7}. (\textit{Banach fixed point theorem}) [6]\textit{. Let}
$T:M\subseteq X\rightarrow M$ \textit{be a contractive operator, where}
$(X,d)$ \textit{is a complete metric space \ and} $M\subseteq X$ \textit{is a
closed nonempty set. Then the following holds}: \textit{(a)} \textit{Existence
and uniqueness. The equation} $x=Tx$, $x\in M$, \textit{has exactly one}
\textit{solution, i.e}., $T$ \textit{has exactly a fixed point} $x$
\textit{on} $M$. \textit{(b) Convergence of the iteration. For any fixed}
$x_{0}\in M$, \textit{the sequence} $\left(  x_{n}\right)  _{n\geq0}$
\textit{of} \textit{succesive approximations, where} $x_{n+1}=Tx_{n}$,
$n=0,1,2,...$, \textit{converges to the fixed point} $x$.\medskip

\textbf{Theorem 8}. (Corollary of the Banach fixed point theorem) [6].
\textit{Let }$X$ \textit{be a complete metric space}, $M\subseteq X$
\ \textit{be a closed nonempty set of }$X$, \textit{and let} $P$ \textit{be a
metric space, called the parameter space. Suppose the following conditions are
satisfied: (i) For each }$p\in P$, \textit{the operator} $T_{p}:M\subseteq
X\rightarrow M$ \textit{is }$k-$\textit{contractive},\textit{ but with }%
$\ $\textit{the same} \textit{contractive constant} $k$\textit{ for}
\textit{all the operators} $T_{p}$. \textit{(ii) for a fixed} $p_{0}\in P$
\textit{and for all }$x\in M$, $\lim_{p\rightarrow p_{0}}T_{p}x=T_{p_{0}}x$.
\textit{Then: (a) For each} $p\in P$, \textit{the equation }$x=T_{p}x$
\textit{has} \textit{exactly one solution} $x_{p}\in M$. (b) $\lim
_{p\rightarrow p_{0}}x_{p}=x_{p_{0}}$.\medskip

\textbf{Theorem 9}. (\textit{Lax-Milgram Theorem}) [8]. \textit{Let}
$B:X\rightarrow X$ \textit{be a linear continuous operator on the Hilbert
space} $X$. \textit{Suppose there is} $c>0$ \textit{such that} $\left\vert
\left(  Bu|u\right)  \right\vert \geq c\left\Vert u\right\Vert ^{2}$
\textit{for all} $u\in X$. \textit{Then, for each given} $f\in X$, \textit{the
operator equation}\quad\ $Bu=f,$\quad$u\in X$ \textit{has a unique solution}.

\begin{center}
\bigskip

\textbf{2. PERTURBED\ OPERATORS\ ON\ BANACH\ SPACES}
\end{center}

Next, we establish our main result.\medskip

\textbf{Lemma 1. } \textit{Let} $\ X$ \textit{be a Banach space over
}$\mathbb{K}=\mathbb{R}$\textit{ or} $\mathbb{C}$\textit{.}

\textit{Suppose that the following conditions are verified}:

\qquad(i) \textit{The operato}r $K:X\rightarrow X$ \textit{is} $k-$%
\textit{contractive} \textit{with} $k\in\lbrack0,1)$.

\qquad(ii) \textit{The operator }$I:X\rightarrow X$\textit{ is the identity.}

\qquad\textit{Then:}

\qquad(a) \textit{The operator }$F:=I+K$ \textit{is a homeomorphism.}

\textit{\medskip}

\textbf{Proof. \ Ad\ (a). }Equation
\begin{equation}
F(x)=y\text{ \quad for any fixed\quad}y\in X\text{,}%
\end{equation}
is equivalent to equation,
\[
x=y-K(x).
\]
Define, the operator
\[
T_{y}:X\rightarrow X,\quad T_{y}(x):=y-K(x).
\]
Operators $\left\{  T_{y}\right\}  _{y\in X}$ fulfill the hypotheses of
Theorem 8. In fact: $X$ is the closed nonempty set of the complete metric
space $X$. The parameter space is also $X$. Furthermore (i) \ For fixed
parameter $y\in X$, operator $T_{y}$ is $k-$contractive and the contraction
constant $k$ is the same for all operators $\left\{  T_{y}\right\}  _{y\in X}%
$, which is the contraction constant of operator $K$. (ii) For any fixed
parameter $y_{0}\in X$ \ and for all $x\in X$
\[
\lim_{y\rightarrow y_{0}}T_{y}(x)=T_{y_{0}}(x)\text{.}%
\]
Then the conclusion of Theorem 8 are verified, i.e.:

Conclusion (a). For each fixed $y$ in the parameter space $X$, there exists a
unique $x_{y}$ \ in the Banach space $X$, such that $T_{y}(x_{y})=x_{y}$, or
written in other way%
\[
y-K(x_{y})=x_{y}\text{.}%
\]
Hence for any fixed $y\in X$, there is a unique $\ x_{y}\in$\ $X$ wich is
solution of equation (1), i.e.,
\[
F(x_{y})=y.
\]
Hence $F$ is a bijection, where $F^{-1}(y)=x_{y}$.

Conclusion (b). For any fixed parameter $y_{0}\in X$, we have
\[
\lim_{y\rightarrow y_{0}}x_{y}=x_{y_{0}}\Rightarrow\lim_{y\rightarrow y_{0}%
}F^{-1}(y)=F^{-1}(y_{0}).
\]
Therefore the operator $F^{-1}$ is continuous.

Since the operator $K$ is continuous, being $k-$contractive, and since
operator $I$ is continuous, then the operator $\ F$ is continuous. Since $F$
is a bijection, and since $F$ and $F^{-1}$ are continuous operators, therefore
$F$ is a homeomorhism.

\medskip

\textbf{Lemma 2. } \textit{Let} $\ X$ \textit{be a Banach space over
}$\mathbb{K}=\mathbb{R}$\textit{ or} $\mathbb{C}$\textit{, and let \ }

\textit{\ }$f:X\rightarrow X$, $f:=F+C$ \textit{be a weakly coercive operator.
Suppose that:}

\qquad(i) $F:X\rightarrow X$\textit{ is a proper continuous opertor.}

\qquad(ii) $C:X\rightarrow X$\textit{ is a compact operator.}

\qquad\textit{Therefore:}

\qquad(a) \textit{Operator} $f$ \textit{is proper}.

\medskip

\textbf{Proof. Ad\ (a). }Since $X$ is a Banach space, a set in $X$ is compact
if, and only if, it is sequentially compact. In the following we do not
distinguish the notation between sequences and their subsequences. Let $A$ be
any fixed compact \ set of $X$; so $A$ is a bounded and closed set of $X$. Let
fix any sequence $(y_{n})_{n\geq1}$ in $A$,
\[
(y_{n})_{n\geq1}\subset A\text{,\quad where\quad}y_{n}=f(x_{n})\text{ for all
}n.
\]
Since $A$ is a bounded set, and since $f$, by hypothesis, is a weakly coercive
operator, then $f^{-1}(A)$ is a bounded set and the sequence
\[
(x_{n})_{n\geq1}\subset f^{-1}(A)
\]
is a bounded set.

Since $A$ is a compact set and $(f(x_{n}))_{n\geq1}\subset A$, there is a
subsequence $(x_{n})_{n\geq1}$ such that \ the subsequence $(f(x_{n}%
))_{n\geq1}$ converges to a point $u$\ of $A$, i.e.,
\[
\lim_{n\rightarrow\infty}f(x_{n})=u\in A\text{.}%
\]
\ 

Since the subsequence $(x_{n})_{n\geq1}$ is a bounded set and $C$ is a compact
operator, there is a subsequence $(x_{n})_{n\geq1}$ such that $(C(x_{n}%
))_{n\geq1}$ converges to a point $v$ of $X$, i.e.,
\[
\lim_{n\rightarrow\infty}C(x_{n})=v\in X\text{.}%
\]
Then there is a subsequence $(x_{n})_{n\geq1}$ such that%
\[
\lim_{n\rightarrow\infty}(f-C)(x_{n})=u-v\in X
\]
or equivalently%
\[
\lim_{n\rightarrow\infty}F(x_{n})=u-v.
\]

Since $F$ is a proper operator and since \ $\left\{  (F(x_{n})_{n\geq1}%
)\cup\left\{  u-v\right\}  \right\}  $ is a compact set, then
\[
F^{-1}\left\{  (F(x_{n}))_{n\geq1}\cup\left\{  u-v\right\}  \right\}
\]
is a compact set. Hence there is a convergent subsequence of the set%
\[
\left\{  \left(  x_{n}\right)  _{n\geq1}\cup F^{-1}(u-v)\right\}
\]
that we call again $(x_{n})_{n\geq1}$, i.e.,
\[
\lim_{n\rightarrow\infty}x_{n}:=w\text{.}%
\]
Furthermore, since the operator $f$ is continuous as addition of continuous
operators, and since the set $\ A$ is closed, therefore $f^{-1}(A)$ is also a
closed set. Now
\[
\lim_{n\rightarrow\infty}x_{n}=w\text{ \qquad and\qquad\ }(x_{n})_{n\geq
1}\subset f^{-1}(A),
\]
then $w\in$ $f^{-1}(A)$. Therefore the sequence
\[
(x_{n})_{n\geq1}\subset f^{-1}(A)
\]
has a subsequence $(x_{n})_{n\geq1}$ which converges to a point $w\in
f^{-1}(A)$, and hence $f^{-1}(A)$ is a compact set. If no points in the
compact set $A$ \ belong to $f(X)$, then $f^{-1}(A)=\emptyset$, which is a
compact set. If there is only a finite number of points $f(x_{i}),i=1,..,N$ in
$A$, trivially, its preimage by \ $f$, $x_{i},i=1,...,N$, is a compact set.
Therefore the preimage of any compact set $A$ for $f$ is also a compact set.
Hence $\ f$ is a proper operator.\medskip

\textbf{Lemma 3. }\textit{Let} $\ X$ \textit{be a Banach space over
}$\mathbb{K}=\mathbb{R}$\textit{ or} $\mathbb{C}$\textit{, and let be}

\textit{\ }$f:X\rightarrow X$, $f:=F+C$\textit{. Suppose that:}

\qquad(i) $F:X\rightarrow X$ \textit{is a} $C^{1}-$\textit{operator, and there
is} $x_{0}\in X$ \textit{such that},

\qquad\ $F^{\prime}(x_{0})$ \textit{is \ bijective}.

\qquad(ii) $C:X\rightarrow X$\textit{ is a} $C^{1}-$\textit{compact operator}

\qquad\textit{Then:}

\qquad(a) $f$ \textit{is a} $C^{1}-$\textit{Fredholm operator of index zero,
i.e.},

\qquad$f^{\prime}(x)\in\mathcal{F}(X,X)$, \textit{and} $\operatorname*{Ind}%
(f^{\prime}(x)=0$ \textit{for all} $x\in X$.

\bigskip

\textbf{Proof.} \textbf{Ad\ (a). }Since the linear continuos operator
$F^{\prime}(x_{0})$ is a bijection%

\[
\ker(F^{\prime}(x_{0}))=\left\{  0\right\}  \text{,\quad and\quad
}\operatorname{R}(F^{\prime}(x_{0}))=X\text{;}%
\]
and so
\[
\dim(\ker(F^{\prime}(x_{0})))=0\text{, and }\operatorname{codim}%
(\operatorname{R}(F^{\prime}(x)))=0.
\]
Therefore, $\ $%
\[
F^{\prime}(x_{0})\in\mathcal{F}(X,Y)\text{, and }\operatorname*{Ind}%
(F^{\prime}(x_{0}))=0.
\]
Since $C$ is a $C^{1}-$compact operator, Theorem 4 implies that $C^{\prime
}(x)$ is also a compact operator. Since $F^{\prime}(x_{0})\in\mathcal{F}%
(X,Y)$, with $\operatorname*{Ind}(F%
\acute{}%
(x_{0}))=0$ and since $C^{\prime}(x_{0})\in\mathcal{L}(X,Y)$ is a compact
operator, Theorem 5 implies that \ the operator
\begin{equation}
f^{\prime}(x_{0})=F^{\prime}(x_{0})+C^{\prime}(x_{0})\in\mathcal{F}(X,Y)\text{
and \ }\operatorname*{Ind}(f^{\prime}(x_{0}))=0\text{.}%
\end{equation}
Since $f$ is a $C^{1}-$operator as addition of $C^{1}-$operators, then the
operator%
\[
f^{\prime}:X\rightarrow\mathcal{L}(X,X)\text{, }x\longmapsto f^{\prime
}(x)\text{,}%
\]
is continuous\textit{. }Since the mapping \ $x\longmapsto f^{\prime}(x)$ is
continuous, Theorem 6 implies that, the index $\operatorname*{Ind}(f^{\prime
}(x))$ is locally constant. Hence, since domain of $f$ is $\ X$, the index is
independent of \ $x\in X$. \ Therefore, equation (2) implies that%
\[
\operatorname*{Ind}(f)=\operatorname*{Ind}(f^{\prime}(x_{0}))=0\text{.}%
\]
Hence $f$ is a $C^{1}-$Fredholm operator of index zero.

\medskip

\textbf{Theorem 10. }\textit{Let} $\ X$ \textit{be a Banach space over
}$\mathbb{K}=\mathbb{R}$\textit{ or} $\mathbb{C}$\textit{, and let \ }

\textit{\ }$f:X\rightarrow X$, $f:=F+C$ \textit{be a weakly coercive operator.
Suppose that:}

\qquad(i) $F$\textit{ is a} $C^{1}-$\textit{proper operator, and there is}
$x_{0}\in X$ \textit{such that}

\qquad\ \ \ \ $F^{\prime}(x_{0})$ \textit{is a linear homeomorphism.}

\qquad(ii) $C$\textit{ is a }$C^{1}-$\textit{compact \ operator}.

\qquad\textit{Therefore}:

\qquad(a) \textit{If for any fixed }$x\in X$, $f^{\prime}(x)$\textit{ is
weakly coercive, then} $f$ \textit{is a global}

\qquad$C^{1}-$\textit{diffeomorphism.}

\qquad(b)\textit{ If }$y\in X$ \textit{is a regular value of} $f$,\textit{
then the equation} $f(x)=y$\textit{ has }

\qquad\textit{at most finitely many solutions}

\bigskip

\textbf{Prooof. }Since $F$ is a proper continuous operator, $C$ a compact
operator, and $f$ a weakly coercive operator, then Lemma 2 implies that $f$
\ is a proper operator.

Hypothesis (i), (ii), and Lemma 3 imply that is a $C^{1}-$Fredholm operator of
index zero, i.e.%
\begin{equation}
f^{\prime}(x)\in\mathcal{F}(X,X)\text{, and }\operatorname*{Ind}(f^{\prime
}(x))=0\text{ for any fixed }x\in X\text{.}%
\end{equation}

\textbf{Ad\ (a).} The statement that for any fixed $x\in X$, the operator
$f^{\prime}(x)$ is weakly coercive, will drive us to
\[
\ker(f^{\prime}(x))=0\text{ for all }x\text{ in }X\text{.}%
\]
In fact if $\ \ker(f^{\prime}(x))\neq\left\{  0\right\}  $, there would be an
$y\in X$ with $y\neq0$ such that
\[
\left(  f^{\prime}(x)\right)  (y)=0\text{.}%
\]
Let $\lambda\in\left(  0,+\infty\right)  $, then%
\[
\left\Vert \lambda y\right\Vert =\left\vert \lambda\right\vert \left\Vert
y\right\Vert \rightarrow\infty\text{\quad when\quad}\lambda\rightarrow
\infty\text{,}%
\]
and
\[
\left\Vert \left(  f^{\prime}(x)\right)  \left(  \lambda y\right)  \right\Vert
=\left\vert \lambda\right\vert \left\Vert \left(  f^{\prime}(x)\right)
\left(  y\right)  \right\Vert =0\text{\quad when\quad}\lambda\rightarrow
\infty\text{,}%
\]
since $\left\vert \lambda\right\vert \left\Vert f^{\prime}(x)y\right\Vert =0$
for all $\lambda$, because of $y\in\ker(f^{\prime}(x))\backslash\left\{
0\right\}  $. This is a contradiction, because $f^{\prime}(x)$ is weakly
coercive. Therefore
\begin{equation}
\ker(f^{\prime}(x))=\left\{  0\right\}  \quad\text{for all }x\in X\text{.}%
\end{equation}
Hence $f^{\prime}(x)$ is an injective operator for any fixed $x\in X.$

Equations (3) and (4) imply that%
\[
\operatorname{codim}R(f^{\prime}(X))=X\text{, and then }\operatorname{R}%
(f^{\prime}(X))=X.
\]
Hence $f^{\prime}(x)$ is a bijective linear continuous operator for any fixed
$x\in X$. Theorem 2 implies that $f$ is a local $C^{1}-$diffeomorphism at any
$x\in X$.

Due to $f$ \ being a proper operator and a local $C^{1}-$diffeomorphism at any
$x\in X$, Theorem 1 implies that $f$ is a global $C^{1}-$diffeomorphism.

\textbf{Ad (b). }Since $y$ is a regular value for $f$, then $f^{-1}(y):=G$ is
empty or consists solely of regular points $\ G:=\left\{  x_{i}\right\}  $.
Then if $G\neq\emptyset$, there is $x_{i}$, such that $f$ has a submersion at
$x_{i}$, that is: $f(x_{i})=y$, $\ker(f^{\prime}(x_{i}))$ has topological
complement in $X$ and $f^{\prime}(x_{i})$ is surjective. Since $f^{\prime
}(x_{i})$ is a surjective operator, $\operatorname{R}(f^{\prime}(x_{i}))=X$,
and then $\operatorname{codim}\operatorname{R}(f^{\prime}(x_{i}))=0$. Equation
(3) implies that $f^{\prime}(x)\in\mathcal{F}(X,X)$, and $\operatorname*{Ind}%
(f^{\prime}(x)=0$ for any fixed $x\in X$. Since $\operatorname*{Ind}%
(f^{\prime}(x_{i}))=0$ and since $\operatorname{codim}(\operatorname{R}%
(f^{\prime}(x_{i})))=0$, therefore $\dim($ $\ker(f^{\prime}(x_{i})))=0$. Hence
$\ker(f^{\prime}(x_{i}))=\left\{  0\right\}  $ and $f^{\prime}(x_{i})$ is
injective. Thus $f^{\prime}(x_{i})$ is a bijection. Since $f^{\prime}(x_{i})$
is a bijection, then the local inverse mapping theorem 2 \ implies that $f$ is
a local $C^{1}-$diffeomorphism at $x_{i}$. Therefore there is a neighborhood
$U(x_{i})$ of $x_{i}$, and a neighborhood $V(y)$ of $y$, such that, $f(x)=y$
with $x\in U(x_{i})$ if and only if $x=x_{i}$.

Since $f$ is a proper operator and $y$ is a compact set, then $G=$ $f^{-1}(y)$
is a compact set. If $\ G\neq\emptyset$,\ since there is a finite subcovering
$\cup_{1\leq i\leq N}U(x_{i})$\ of $G$ for the open covering $\cup_{x_{i}\in
G}U(x_{i})$ of the compact set $G$, therefore there is at most a finite number
of points $x$ such that $f(x)=y.\medskip$

\textbf{Corollary 1.} \textbf{ }\textit{Let} $\ X$ \textit{be a Banach space
over }$\mathbb{K}=\mathbb{R}$\textit{ or} $\mathbb{C}$\textit{, and let \ }

\textit{\ }$f:X\rightarrow X$, $f:=F+C$, \textit{ where}:

\qquad(i) $F$ \textit{is a linear homeomorphism}

\qquad(ii) $C$\textit{ is a linear compact operator.}

\qquad(iii) $\left\Vert F\left\Vert \neq\right\Vert C\right\Vert $

\qquad\textit{Then:}

\qquad(a) \ $f$ \textit{is a global} $C^{1}-$\textit{diffeomorphism.}\bigskip

\textbf{Prooof. Ad\ (a). }Since $f$, $F$, $C$ are linear continuous operators,
then they are $C^{1}-$operators, and for any fixed $x\in X$%
\[
\left(  F^{\prime}(x)\right)  (y)=F(y)\text{, }\left(  C^{\prime}(x)\right)
(y)=C(y)\text{, and }\left(  f^{\prime}(x)\right)  (y)=f(y)\text{ for all
}y\in X\text{.}%
\]
Therefore%
\[
\left\vert \left\Vert F(\frac{y}{\left\Vert y\right\Vert })\right\Vert
-\left\Vert C(\frac{y}{\left\Vert y\right\Vert })\right\Vert \right\vert
\left\Vert y\right\Vert =\left\vert \left\Vert F(y)\right\Vert -\left\Vert
C(y)\right\Vert \right\vert \leq\left\Vert f(y)\right\Vert
\]%
\[
\leq\left(  \left\Vert F(\frac{y}{\left\Vert y\right\Vert })\right\Vert
+\left\Vert C(\frac{y}{\left\Vert y\right\Vert })\right\Vert \right)
\left\Vert y\right\Vert \quad\text{for all }y\neq0.
\]
Hence%
\[
\left\vert \left\Vert F\right\Vert -\left\Vert C\right\Vert \right\vert
\left\Vert y\right\Vert \leq\left\Vert f(y)\right\Vert \leq\left(  \left\Vert
F\right\Vert +\left\Vert C\right\Vert \right)  \left\Vert y\right\Vert ,
\]
and so if $\left\Vert F\right\Vert \neq\left\Vert C\right\Vert $,%
\[
\left\Vert f(y)\right\Vert \rightarrow\infty\quad\text{if and only if \quad
}\left\Vert y\right\Vert \rightarrow\infty.
\]
Thus $f$ is weakly coercive, and furthermore $f^{\prime}(x)$ is also weakly
coercive for any fixed $x\in X$.

Since $F$ is a homeomorphism and since continuous operators leave invariant
compact sets, $F$ is a proper operator.

Since $F$ is a linear homeomorphism, $F^{\prime}(x)=F$ and $\left(
F^{-1}\right)  ^{\prime}(x)=F^{-1}$ for all $x\in X$. Therefore $F^{\prime
}(x)$ is a linear homeomorphism for any fixed $x\in X$.

Hipotheses of Theorem 10 (a) are fulfilled, then $f$ is a global $C^{1}-$diffeomorphism.

\bigskip

\textbf{Corollary 2. }\textit{Let} $\ X$ \textit{be a Banach space over
}$\mathbb{K}=\mathbb{R}$ or $\mathbb{C}$,

\textit{and let }$f:X\rightarrow X$, \ $f=I+K+C$ \textit{be} \textit{a weakly
coercive operator},

\qquad\textit{where}:

\qquad(i) $I$ \textit{is the identity operator.}

$\qquad$(ii) $K$ \textit{is a linear} $k-$\textit{contractive operator},
\textit{with} $k\in\lbrack0,1)$.

$\qquad$(iii) $C$ \textit{is a }$C^{1}-$\textit{compact operator.}

\qquad(iv) $I+K+C^{\prime}(x)$ \textit{is a weakly coercive operator for all}
$x\in X$.

\qquad\textit{Therefore:}

\qquad(a) $f$ \textit{is a global} $C^{1}-$\textit{diffeomorphism}.

\qquad(b) \textit{If operator }$K_{1}=-K$, \textit{and operator} $C_{1}=-C$,
\textit{then} $K_{1}+C_{1}$ \textit{has}

\qquad\qquad\textit{a unique fixed point.\medskip}

\textbf{Proof. Ad (a).\ }Set the operator
\[
F:X\rightarrow X,\quad F(x):=I(x)+K(x).
\]
Therefore the weakly coercive operator $f$ \ can be written as
\[
f=F+C.
\]

Lemma 1 implies that $F$ is a homeomorphism, and then a proper operator,
because of $F^{-1}$ being continuous, and since cotinuous operators map
compact sets into compact sets.

Since $I$ and $K$ are linear, the homeomorphism $F$ is linear. Therefore
$F^{\prime}(x)=F$ and $\left(  F^{-1}\right)  ^{\prime}(x)=F^{-1}$ for all
$x\in X$. Hence $F^{\prime}(x)$ is a linear homeomorphism for any fixed $x\in
X.$

The hypothesis (iv) implies that $f^{\prime}(x)$ is a weakly coercive operator
for any fixed $x\in X$.

Since hypotheses of Theorem 10 (a) are fulfilled, $f$ is a global $C^{1}%
-$diffeomorphism$.$

\textbf{Ad (b). }Since $f$ \ is a global diffeomorfism, there is a unique
$x\in X$ with $f(x)=0$, therefore $\left(  K_{1}+C_{1}\right)  (x)=x$%
.$\medskip$

\textbf{Corollary 3. }\textit{Let} $\ X$ \textit{be a Hilbert space space over
}$\mathbb{K}=\mathbb{R}$ or $\mathbb{C}$,

\textit{and let }$f:X\rightarrow X$, \ $f=F+C$ \textit{be} \textit{a weakly
coercive operator.}

\qquad\textit{Suppose that}:

\qquad(i) $F$\textit{ is a linear continuous operator \ and there is a
constant} $k>0$

\qquad such that \ $\left\vert \left(  Fx\mid x\right)  \right\vert \geq
k\left\Vert x\right\Vert ^{2}$ for all $x\in X.$

\qquad(ii) $C$\textit{ is a }$C^{1}-$\textit{compact operator, for all }$x\in
X$, $\left\Vert F\left\Vert \neq\right\Vert C^{\prime}(x)\right\Vert $.

\qquad Therefore:

\qquad(a) $f$ \textit{is a global} $C^{1}-$\textit{diffeomorphism.}\medskip

\textbf{Proof. Ad (a).\ } Theorem 9 implies that $F$ is a linear continuous
bijection. Then the open mapping Theorem 3 (ii)\ implies that the linear
operator $F^{-1}$ is also continuous. Therefore $F$ is a linear homeomorphism,
and then a linear $C^{1}-$proper operator.

Since $F$ is linear and continuous, $(F^{\prime}(x))(y)=F(y)$ for any fixed
$x\in X$ and all $y\in X$. \ Therefore for any fixed $x\in X$, $F^{\prime}(x)$
is a linear homeomorphism.

For any fixed $x\in X$, \ operator $f^{\prime}(x)$ is weakly coercive. In fact:%

\[
\left(  \left\Vert F\left(  \frac{y}{\left\Vert y\right\Vert }\right)
\right\Vert +\left\Vert \left(  C^{\prime}(x)\right)  \left(  \frac
{y}{\left\Vert y\right\Vert }\right)  \right\Vert \right)  \left\Vert
y\right\Vert \geq\left\Vert \left(  f^{\prime}(x)\right)  \left(  y\right)
\right\Vert =\left\Vert F(y)+\left(  C^{\prime}(x)\right)  \left(  y\right)
\right\Vert
\]%
\[
=\left\Vert F\left(  \frac{y}{\left\Vert y\right\Vert }\right)  \left\Vert
y\right\Vert +\left(  C^{\prime}(x)\right)  \left(  \frac{y}{\left\Vert
y\right\Vert }\right)  \left\Vert y\right\Vert \right\Vert \geq\left\vert
\left\Vert F\left(  \frac{y}{\left\Vert y\right\Vert }\right)  \right\Vert
-\left\Vert \left(  C^{\prime}(x)\right)  \left(  \frac{y}{\left\Vert
y\right\Vert }\right)  \right\Vert \right\vert \left\Vert y\right\Vert
\text{.}%
\]
for all $y\in Y$. Then%
\[
\left(  \left\Vert F\right\Vert +\left\Vert C^{\prime}(x)\right\Vert \right)
\left\Vert y\right\Vert \geq\left\Vert \left(  f^{\prime}(x)\right)  \left(
y\right)  \right\Vert \geq\left\vert \left\Vert F\right\Vert -\left\Vert
C^{\prime}(x)\right\Vert \right\vert \left\Vert y\right\Vert \text{.}%
\]
Therefore
\[
\left\Vert \left(  f^{\prime}(x)\right)  \left(  y\right)  \right\Vert
\rightarrow\infty\quad\text{if and only if }y\rightarrow\infty\text{.}%
\]
Hence for any fixed $x\in X$, $f^{\prime}(x)$ is weakly coercive .

Theorem 10 (a) implies that $f$ is a global $C^{1}-$diffeomorphism.

\bigskip

\begin{center}
\textbf{3}. \textbf{EXAMPLES.\medskip}
\end{center}

Let $G$ be a nonempty bounded open set in $\mathbb{R}^{n}$, and $\Omega
:=\overline{G}$. In the two following examples let \ $X$ be the Banach space
$X=C\left(  \Omega\right)  $ of continuous functions on the compact set
$\Omega$ with the norm $\parallel u\parallel=\max_{x\in\Omega}\left\vert
u(x)\right\vert $, and let the operator $f$ $:X\longrightarrow X$ \ be the
perturbed operator $f:=F+C$, where operators $F,C:X\rightarrow X$ will be
defined below. We want to know if $f$ \ is a $C^{1}-$diffeomorphism and then
if for any $v\in X$ \ there is one and only one $u\in X$ such that $fu=v$.
\medskip

\textbf{Example 1. }Let be \textbf{\textit{ }}$F:=I+K$, with $I$ the identity
operator on the space $X$, and \ $K,C:X\rightarrow X$ a $k-$contractive and a
compact operator respectively to be defined below.\textbf{ }We define these
particular operators $K,C$ to be consider.

Operator $K$ is defined in the following way:
\[
\text{for any }u\in X\text{,\quad}u\mapsto Ku\text{,\quad where\quad
}Ku(x):=\int_{\Omega}k(x,y)u(y)dy\text{ for all }x\text{ in }\Omega\text{.}%
\]
being
\[
k(\cdot,\cdot):\Omega\times\Omega\rightarrow\mathbb{R}%
\]
a continuous function with $\ \max_{\left(  x,y\right)  \in\Omega\times\Omega
}\left\vert k(x,y)\right\vert \operatorname*{meas}\left(  \Omega\right)  <1$.

Since $\Omega\times\Omega\subset\mathbb{R}^{n}\times\mathbb{R}^{n}$ is a
compact set \ and $k(\cdot,\cdot)$ a continuous function, there is
$\max_{\left(  x,y\right)  \in\Omega\times\Omega}\left\vert k(x,y)\right\vert
:=M$, and furthermore $k(\cdot,\cdot)$ is uniformly continuous on
$\Omega\times\Omega$. Trivially $K$ maps $X$ into $X$, and it is a linear operator.

Since for $u\in X$,%
\[
\left\Vert Ku\right\Vert =\max_{x\in\Omega}\left\vert
{\displaystyle\int_{\Omega}}
k(x,y)u(y)dy\right\vert \leq\max_{x\in\Omega}%
{\displaystyle\int_{\Omega}}
\left\vert k(x,y)u(y)\right\vert dy\leq M\operatorname*{meas}\left(
\Omega\right)  \left\Vert u\right\Vert \text{,}%
\]
then the operator $K$ is bounded. Since $K$ is a linear and bounded operator,
it is a $C^{1}-$operator with $K%
\acute{}%
(u)=K$ for all $u\in X$.

$K$ is a contractive operator since $M\operatorname*{meas}\left(
\Omega\right)  :=k<1$. In fact, for
\[
u,v\in X\text{,}\left\Vert Ku-Kv\right\Vert =\left\Vert K(u-v)\right\Vert
=\max_{x\in\Omega}\left\vert \int_{\Omega}(k(x,y)(u-v)(y))dy\right\vert
\]%
\[
\leq M\operatorname*{meas}\left(  \Omega\right)  \left\Vert u-v\right\Vert
=k\left\Vert u-v\right\Vert \text{.}%
\]

Operator $C$ \ is defined in the following way:
\[
\text{for any }u\in X\text{,\quad}u\mapsto Cu\text{,\quad where\quad
}Cu(x):=\int_{\Omega}g(x,y)u(y)dy\text{ for all }x\text{ in }\Omega\text{.}%
\]
being
\[
g(\cdot,\cdot):\Omega\times\Omega\rightarrow\mathbb{R}\text{.}%
\]
a continuous function. The functions $\ k(\cdot,\cdot)$, $g(\cdot,\cdot)$ can
be selected when linear operator $K$ $+$ $C$ verify the condition $\left\Vert
K+C\right\Vert \neq1$, which is required so that the operator $I+K+C$ to be
weakly coercive.

By using the same reasoning that in the previous operator, we obtain that
$C:X\rightarrow X$ is a linear $\ C^{1}(X)$ operator with $C%
\acute{}%
(u)=C$ for all $u$ in $X$, and that there is a real number $N$ such that
$\max_{\left(  x,y\right)  \in\Omega\times\Omega}\left\vert g(x,y)\right\vert
\operatorname*{meas}\left(  \Omega\right)  :=N$, and furthermore
$g(\cdot,\cdot)$ is uniformly continuous on $\Omega\times\Omega$. \ 

Since $\ C$ is a continuous operator, to show that $C$ is a compact operator,
we have to prove that $C$ maps bounded sets into relatively compact sets. To
this end, let $L$ be the bounded set $L:=\left\{  u\in X:\left\Vert
u\right\Vert \leq r\right\}  $. Since $\left\Vert Cu\right\Vert \leq
N\operatorname*{meas}(\Omega)r$ for all $u\in L$, the set $C(L)$ is a bounded
set. Then, by the Arzel\'{a}-Ascoli theorem [6-7], to show that $C(L)$ is a
relatively compact set, we only have to prove that $C(L)$ is equicontinuous.
Let $Cu\in C(L)$ and $x,z\in\Omega$, we have%
\[
\left\vert Cu(x)-Cu(z)\right\vert =\left\vert \int_{\Omega}g(x,y)u(y)dy-\int
_{\Omega}g(z,y)u(y)dy\right\vert
\]%
\[
=\left\vert \int_{\Omega}u(y)(g(x,y)-g(z,y))dy\right\vert \leq r\int_{\Omega
}\left\vert g(x,y)-g(z,y)\right\vert dy\text{.}%
\]
Since $g(\cdot,\cdot)$ is uniformly continuous,for any $\varepsilon>0$ there
is $\delta(\varepsilon)>0$ such that
\[
\left\Vert x-z\right\Vert <\delta(\varepsilon)\Rightarrow\left\vert
g(x,y)-g(z,y)\right\vert <\varepsilon\text{,}%
\]
and then%
\[
\left\vert Cu(x)-Cu(z)\right\vert \leq r\varepsilon\operatorname*{meas}\left(
\Omega\right)  \text{\quad for all }u\in L\text{,}%
\]%
\[
\text{ and all }x,z\in\Omega\text{ with }\left\Vert x-z\right\Vert
<\delta(\varepsilon)\text{.}%
\]
Therefore $C(L)$ is equicontinuous. Hence $C$ is a compact operator.

By Lemma 1 the operator $F$ is a homeomorphism. Trivially $F$ is linear as
addition of linear operators. Therefore $F$ is a $C^{1}-$linear homeomorphism,
and
\[
F^{\prime}(u)=F\text{ \quad for all\quad\ }u\quad\text{ in\quad}X.
\]
Hence operator $F^{\prime}(u)$ is a linear homeomorphism for all $u\in X$.

Operator $f=F+C=I+(K+C):X\rightarrow X$ is a linear continuous operator as
addition of linear continuous operators, then%
\begin{equation}
f^{\prime}(u)=f\text{ \quad for all\quad\ }u\quad\text{ in\quad}X.
\end{equation}
\
\[
\left\Vert fu\right\Vert \leq\left\Vert f\right\Vert \left\Vert u\right\Vert
\text{,}%
\]
$_{{}}$which implies $\left\Vert fu\right\Vert $ is bounded if $\left\Vert
u\right\Vert $ is bounded.

Since $\left\Vert K+C\right\Vert \neq1$,%
\[
\left\Vert f(u)\right\Vert =\left\Vert (K+C+I)u\right\Vert \geq\left\vert
\left\Vert (K+C)\frac{(u)}{\left\Vert u\right\Vert }\right\Vert \left\Vert
u\right\Vert -\left\Vert u\right\Vert \right\vert \text{, for all }%
u\neq0\text{,}%
\]
then,%
\[
\left\Vert f(u)\right\Vert \geq\left\vert \left\Vert K+C\right\Vert
-1\right\vert \left\Vert u\right\Vert
\]
Therefore $\left\Vert f(u)\right\Vert \rightarrow\infty$ as $\left\Vert
u\right\Vert \rightarrow\infty$. Hence $\ f$ is weakly coercive $\ $Formula
(5) implies that \ for any fixed $u\in X$, operator $f^{\prime}(u)$ is also
weakly coercive

Since the hypotheses of Theorem 10 (a) are verified, $f$ is a $C^{1}%
-$diffeomorphim..\medskip

\textbf{Example 2.} In this example the operator $F$ is the identity $I$ on
the space $X$, which is a linear homeomorphism, and then a $C^{1}-$proper
operator with $I^{\prime}(u)=I$ for all $u$ in $X$.

Compact operator $C$ is defined in the following way: Let \ $Q:=\Omega
\times\Omega\times\mathbb{R}:$and let%
\[
h:Q\rightarrow\mathbb{R}%
\]
be a $C^{1}-$function.with $\sup_{(x,y,u)\in Q}\left\vert h(x,y,u)\right\vert
=\mathcal{L}$

Now define operator $C$:$X\rightarrow X$, $u\mapsto Cu$ as $\ Cu(x):=%
{\displaystyle\int_{x\in\Omega}}
h(x,y,u(y))dy$ for all $x\in\Omega$. Indeed if $u\in X$, $Cu\in X$:

Since $h$ is continuous $\forall\varepsilon>0$ there is $\delta(\varepsilon
,u)$ such that $\left\Vert x-z\right\Vert +\left\vert u-v\right\vert
<\delta(\varepsilon,u)\Rightarrow\left\vert h(x,y,u)-h(z,y,v)\right\vert
<\varepsilon$, and then
\[
\left\vert Cu(x)-Cu(z)\right\vert =\left\vert
{\displaystyle\int_{x\in\Omega}}
h(x,y,u(y)-h(z,y,u(y))dy\right\vert \leq\varepsilon\operatorname*{meas}%
(\Omega).
\]
Therefore $Cu\in X$.

We use the technique of linearization of $C(u+m)$ with respect to $m$ to study
the existence of first Fr\'{e}chet derivative at the point $u\in X$. This
means that for fixed $u\in X$, we search for a decomposition of the form
\[
C(u+m)=Cu+dC(u;m)+o(\left\Vert m\right\Vert )\text{\quad for all\quad}m\in X
\]
in a neighborhood of zero, where $dC(u;m)$ represents the linear part with
respect to $m$. This linear part should be continuous.\ Indeed,
\[
C(u+m)(x)=\int_{\Omega}h(x,y,u(y)+m(y))\text{ }dy
\]%
\[
=\int_{\Omega}\left\{  h(x,y,u(y)+h_{u}(x,y,u(y))m(y)+r_{2}\right\}  \text{
}dy
\]%
\[
=Cu(x)+\int_{\Omega}h_{u}(x,y,u(y))m(y)\text{ }dy+R_{2}\text{,\quad where
}\left\Vert R_{2}\right\Vert =o\left(  \left\Vert m\right\Vert \right)
\text{, }m\rightarrow0,
\]
for all \ $m$ in a neighborhood of zero. Therefore the linear part is
\[
dC(u;m)(x)=\int_{\Omega}h_{u}(x,y,u(y))m(y)dy,
\]
which certainly is continuous for fixed $\ u\in X$. Hence the operator $C$ is
Fr\'{e}chet differentiable \ at $u$, and the Fr\'{e}chet derivative of $C$ at
$u$ is $C^{\prime}(u)m=dC(u;m)$. Furthermore
\[
\left\Vert C^{\prime}u(m)\right\Vert \leq\mathcal{L}\operatorname*{meas}%
(\Omega)\left\Vert m\right\Vert ,\quad\text{and}\quad\left\Vert C^{\prime
}(u)\right\Vert \leq\mathcal{L}\operatorname*{meas}(\Omega)\mathcal{\ }%
\text{for \ all \ }u\in X.
\]
Since the Fr\'{e}chet derivative \ $C^{\prime}(u)$ exists for all $u\in X$,
then the Fr\'{e}chet derivative of $C$ at $X$, $C^{\prime}:X\rightarrow
\mathcal{L}(X,X)$, also exists. To show that $C\in C^{1}(X,X)$, we have to
prove that $C^{\prime}$ is continuous at any point $u\in X$. Let $u,v,m\in X$.
Since $C^{\prime}(u)$, $C^{\prime}(v)\in\mathcal{L}(X,X)$, then $C^{\prime
}(u)-C^{\prime}(v)\in\mathcal{L}(X,X)$ and
\[
\left\Vert C^{\prime}(u)-C^{\prime}(v)\right\Vert =\sup\limits_{m\in
X,\left\Vert m\right\Vert \leq1}\left\Vert (C^{\prime}(u)-C^{\prime
}(v))m\right\Vert \text{.}%
\]
Since $h_{u}$ is continuoous, because of being $h\in C^{1}(Q)$, $\forall
\varepsilon>0$ there is $\delta(\varepsilon,u)$ such that $\left\vert
u(y)-v(y)\right\vert <\delta(\varepsilon,u)\Rightarrow\left\vert
h_{u}(x,y,u(y))-h_{u}(x,y,v(y))\right\vert <\varepsilon$. Now, observe that if
$\left\Vert u-v\right\Vert <\delta(\varepsilon,u)$ then%
\[
\left\Vert C^{\prime}(u)m-C^{\prime}(v)m\right\Vert =\max\limits_{x\in\Omega
}\left\vert
{\displaystyle\int_{\Omega}}
\left(  h_{u}(x,y,u(y))-h_{u}(x,y,v(y)\right)  m(y)dy\right\vert
\]%
\[
\leq\varepsilon\operatorname*{meas}(\Omega)\left\Vert m\right\Vert \text{.}%
\]
Hence the operator $C^{\prime}$ is continuous, and thus $C$ is a $C^{1}-$operator.

Operator $C$ is compact: In fact, it is continuous since it is a $C^{1}%
-$operator, and furthermore maps bounded sets into relatively compact sets as
we will show

Let $P$ be any bounded subset of $X$, then there is $R>0$ such that $u\in P$
implies that $u\in X$ and $\left\Vert u\right\Vert \leq R$. Let $Q_{1}$ be the
following compact set%
\[
Q_{1}:=\left\{  x,y,u)\in\Omega\times\Omega\times\mathbb{R}:\left\vert
u\right\vert \leq R\right\}  .
\]
To apply the Ascoly-Arzela theorem [6-7], we have to show that $C(P)$ is
bounded and equicontinuous so that $C(P)$ is a relatively compact set. Since
$h$ is continuous and $Q_{1}$ a compact set, then $h$ is uniformly continuous
on $Q_{1}$ and there is $\mathcal{M=}\max_{(x,y,u)\in Q_{1}}\left\vert
h(x,y,u)\right\vert $. Therefore for any fixed $u\in P$, $\left\Vert
u\right\Vert \leq R\Rightarrow\left\vert u(y)\right\vert <R$ for all
$y\in\Omega$, and
\[
\max\limits_{(x,y,u(y))\in Q_{1}}\left\vert
{\displaystyle\int_{\Omega}}
h(x,y,u(y))dy\right\vert \leq\mathcal{M}\operatorname*{meas}(\Omega)\text{.}%
\]
Hence $C(P)$ is bounded. Since $h$ is uniformly continuous on $Q_{1}$,
\[
\forall\varepsilon>0,\text{ }\exists\text{ }\partial(\varepsilon)>0\text{ such
that }\left\Vert x-z\right\Vert <\partial(\varepsilon),x,z\in\Omega
\Rightarrow\left\vert Cu(x)-Cu(z)\right\vert
\]%
\[
\leq%
{\displaystyle\int_{x\in\Omega}}
\left\vert h(x,y,u(y)-h(z,y,u(y))\right\vert <\varepsilon\operatorname*{meas}%
(\Omega)\text{.for all }u\in P.
\]
Therefore $C(P)$ is equicontinuous $\ $Hence $C$ maps bounded sets in
relatively compact sets. Thus $C$ is a compact operator.

Operator $f=I+C$ is weakly coercive. In fact. Since
\[
\left\Vert f(u)\right\Vert \leq\left\Vert u\right\Vert +\left\Vert
Cu\right\Vert \leq\left\Vert u\right\Vert +\mathcal{M}\operatorname*{meas}%
(\Omega),
\]
$f$ is bounded if $u$ is bounded. Since
\[
\left\Vert f(u)\right\Vert \geq\left\vert \left\Vert Cu\right\Vert -\left\Vert
u\right\Vert \right\vert \rightarrow\infty\text{ if }u\rightarrow\infty.
\]
i.e., $\left\Vert f(u)\right\Vert \rightarrow\infty$ if and only if
$\left\Vert u\right\Vert \rightarrow\infty.$

For any fixed $u\in X$, the operator $C^{\prime}u$ is weakly coercive:

Since $\left\Vert (f^{\prime}(u))v\right\Vert \leq\left(  1+\left\Vert
C^{\prime}(u)\right\Vert \right)  \left\Vert v\right\Vert $, then $\left\Vert
(f^{\prime}(u))v\right\Vert $ is bounded if $\left\Vert u\right\Vert $ is
bounded. Since $\left\Vert (f^{\prime}(u))v\right\Vert \geq\left\vert
1-\left\Vert (C^{\prime}(u))\frac{v}{\left\Vert v\right\Vert }\right\Vert
\right\vert \left\Vert v\right\Vert $, then

$\left\Vert (f^{\prime}(u))v\right\Vert \rightarrow\infty$ if $\left\Vert
v\right\Vert \rightarrow\infty$.

Therefore the hypotheses of Theorem 10 (a) are fulfilled. Hence $f$ is a
global $C^{1}-$diffeomorphism

\end{document}